\documentclass[10pt]{article}

\usepackage{graphicx}
\usepackage{subcaption}

\usepackage[a4paper, left=34mm,right=34mm,top=34mm,bottom=34mm]{geometry}
\usepackage[utf8]{inputenc}
\usepackage[T1]{fontenc}
\usepackage[english]{babel}

\usepackage{enumerate}
\usepackage{graphicx}
\usepackage{hyperref}
\hypersetup{
    colorlinks=true,
    linkcolor=blue,
    filecolor=magenta,      
    urlcolor=cyan,
}

\usepackage{mathtools,amsthm,amssymb,amsfonts}
\usepackage{algorithm}
\usepackage{algpseudocode}
\makeatletter
\def\thm@space@setup{
  \thm@preskip=10pt \thm@postskip=10pt
}
\makeatother
\theoremstyle{plain}

\theoremstyle{plain}

\theoremstyle{definition}

\theoremstyle{definition}

\theoremstyle{remark}

\theoremstyle{remark}

\usepackage{caption} 
\usepackage{booktabs}
\usepackage{listings}
\usepackage{color}
\definecolor{dkgreen}{rgb}{0,0.6,0}
\definecolor{gray}{rgb}{0.5,0.5,0.5}
\definecolor{mauve}{rgb}{0.58,0,0.82}
\usepackage{enumitem}

\newcommand{\abs}[1]{\left\lvert#1\right\rvert}
\newcommand{\norm}[1]{\left\lVert#1\right\rVert}
\newcommand{\overbar}[1]{\mkern 3.0mu\overline{\mkern-3.0mu#1\mkern-1.0mu}\mkern 1.0mu}

\DeclareMathOperator*{\argmin}{arg\,min}
\DeclareMathOperator{\diag}{diag}
\DeclareMathOperator{\sign}{sign}
\DeclareMathOperator{\vspan}{span}

\newcommand{\email}[1]{\protect\href{mailto:#1}{#1}}

\usepackage{xcolor}[2.11]
\colorlet{inlinkcolor}{red!100!black}
\colorlet{exlinkcolor}{blue!100!black}
\if@hidelinks
\hypersetup{ hidelinks = true }
\else
\hypersetup{
  colorlinks = true,
  allcolors = inlinkcolor,
  urlcolor = exlinkcolor,
}
\fi

\newenvironment{@abssec}[1]{
        \vspace{.05in}\parindent .0in
        {\upshape\bfseries #1. }\ignorespaces
    }
    {\par\vspace{.1in}}
\renewenvironment{abstract}{\begin{@abssec}{\abstractname}}{\end{@abssec}}
\newenvironment{keywords}{\begin{@abssec}{Keywords}}{\end{@abssec}}

\usepackage{fancyhdr}

\author{
  {\normalsize Qinmeng Zou}\thanks{School of Science, Beijing University of Posts and Telecommunications, Beijing {\rm 100876}, China
    (\email{zouqinmeng@bupt.edu.cn}).}
}
\title{GMRES algorithms over 35 years}
\date{}

\begin{document}
\maketitle
\thispagestyle{fancy}

\begin{abstract}
This paper is about GMRES algorithms for the solution of nonsingular linear systems.
We first consider basic algorithms and study their convergence.
We then focus on acceleration strategies and parallel algorithms that are useful for solving challenging systems.
We also briefly discuss other problems, such as systems with multiple right-hand sides, shifted systems, and singular systems.
\end{abstract}


\section{Introduction}

Generalized minimal residual (GMRES) algorithms are widely used for solving nonsymmetric linear systems arising from partial differential equations.
They possess some optimality property and behave reasonably well in practice.
The first and most well-known algorithm is that of Saad and Schultz, which was initially introduced in a technical report in 1983, and formally published in 1986~\cite{Saad1986}.
Since then, numerous variants appeared, as well as studies of their convergence and accuracy.
In this paper we concentrate on GMRES algorithms.
We sketch the main developments and focus on the algorithmic innovations in the past 35 years.

Consider the nonsingular linear system
\begin{equation}\label{eq:ls}
Ax=b,
\end{equation}
where $A\in\mathbb{R}^{N\times N}$ and $b\in\mathbb{R}^N$.
Given an initial approximation~$x_0$, a projection method constructs a sequence of approximate solutions~$x_n$, $n=1,2,\dots$, such that
\[
x_n \in x_0 + \mathcal{S}_n,\quad r_n \perp \mathcal{C}_n,
\]
where $r_n=b-Ax_n$ is the residual vector and $\mathcal{S}_n$ and $\mathcal{C}_n$ are $n$-dimensional subspaces, called respectively search space and constraint space.
As $A$ is nonsingular, $x_n$ is uniquely defined.
In general, we build a sequence of nested search spaces
\begin{equation}\label{eq:nested}
\mathcal{S}_1 \subset \mathcal{S}_2 \subset \mathcal{S}_3 \subset \dots
\end{equation}
to ensure finite termination; see,~e.g.,~\cite{Liesen2013}.
A minimal residual method obtains the optimal approximation by minimizing the residual norm over all candidate vectors in the search space
\begin{equation}\label{eq:mr0}
\norm{r_n} = \norm{b-Ax_n} = \min_{x\in x_0+\mathcal{S}_n}\norm{b-Ax},
\end{equation}
where $\norm{\cdot}$ denotes the 2-norm or the corresponding induced matrix norm.
This class of methods can be interpreted as a projection process with~$\mathcal{C}_n=A\mathcal{S}_n$.
In other words, they find an optimal correction in~$\mathcal{S}_n$, such that $r_0-r_n$ is the orthogonal projection of~$r_0$ onto~$A\mathcal{S}_n$; see~\cite{Eiermann2000,Saad2000,Eiermann2001} for further analysis.
The Krylov subspace is the most broadly employed search space and is defined by
\begin{equation}\label{eq:krylov}
\mathcal{K}_n = \mathcal{K}_n(A,r_0) = \vspan\{r_0,Ar_0,\dots,A^{n-1}r_0\}.
\end{equation}
These subspaces are nested in the sense of~\eqref{eq:nested} and can be built up gradually using only matrix-vector multiplication by~$A$.
The residual~$r_n$ lies in the next Krylov subspace~$\mathcal{K}_{n+1}$.
If $d$ denotes the grade of~$r_0$ with respect to~$A$, then $A\mathcal{K}_d\subseteq\mathcal{K}_d$ holds.
Therefore, the sequence of Krylov subspaces will eventually become invariant.
A projection method with~$\mathcal{S}_n=\mathcal{K}_n$ is called a Krylov subspace method.
Practically a Krylov subspace method is terminated as soon as the approximate solution is good enough.
Ways to formulate typical Krylov subspace methods and derivations of their properties have been extensively discussed in the literature; see,~e.g.,~\cite{Joubert1990,Freund1992b,Greenbaum1997,Saad2003,vanderVorst2003,Simoncini2007,Liesen2013,Meurant2020} and references therein.
A comparison of Krylov subspace methods can be found in a recent work~\cite{Ghai2019}; see also~\cite{Liesen2013,Saad2020} for a historical perspective.

Classical GMRES algorithms realize the so-called minimal residual Krylov subspace method, or loosely called GMRES method, which is a projection process satisfying~\eqref{eq:mr0} and taking~\eqref{eq:krylov} as search space.
In this case, the projection process can be rewritten as follows:
\begin{subequations}\label{eq:proj}
\begin{align}
& x_n \in x_0 + \mathcal{K}_n, \label{eq:proj1} \\
& r_n \perp A\mathcal{K}_n, \label{eq:proj2}
\end{align}
\end{subequations}
or equivalently,
\begin{equation}\label{eq:mr}
\norm{r_n} = \norm{b-Ax_n} = \min_{x\in x_0+\mathcal{K}_n}\norm{b-Ax}.
\end{equation}
We shall consider both sequential and parallel GMRES algorithms.
Some variants such as restarted or hybrid algorithms exhibit more complex behavior and may no longer possess the finite termination property.
Other crucial aspects that can dramatically affect the performance of Krylov subspace algorithms include the finite precision effect and communication costs; see~\cite{Carson2020b} for a general insight on the cost of the iterative computation and its role in Krylov subspace methods.
In what follows these topics shall be discussed to some extent.
We do not talk here about general-purpose preconditioning techniques like incomplete factorization and algebraic multigrid, since these topics deserve individual coverage by themselves, for which we refer the reader to~\cite{Benzi2002,Wathen2015}; see also~\cite{Ghai2019} for experiments.
For the similar reason, problem-specific techniques will also not be covered, though they can be very successful when prior knowledge of the source problem is available; see,~e.g.,~\cite{Nigro1998,Gander2015b,Pearson2020}.
We only mention certain techniques that are closely related to GMRES iterations.
Throughout most of the paper we restrict to the nonsingular linear system~\eqref{eq:ls}.
In Section~\ref{sec:ext}, however, we briefly mention some work where this restriction is lifted.

The paper is organized as follows.
We begin, in Section~\ref{sec:bac}, by formulating basic GMRES algorithms, giving some equivalent formulations, and presenting theoretical results in both exact and finite arithmetic.
In Section~\ref{sec:acc} we summarize various strategies for improving classical algorithms.
In Section~\ref{sec:pa} we discuss parallel techniques that can reduce or hide communication overhead.
Sections~\ref{sec:ext} is devoted to a brief literature review of related topics.
Finally, concluding remarks are draw in Section~\ref{sec:con}.

\section{Basic algorithms and convergence}\label{sec:bac}

The techniques introduced in this section form the basis for subsequent discussions.
In particular, the algorithm developed by Saad and Schultz~\cite{Saad1986} is still in common use today.
Methodologies have been developed to partially characterize their convergence behavior.
Note that the literature on Krylov subspace methods uses the term ``GMRES'' almost as a synonym for the specific algorithm developed in~\cite{Saad1986}.
For the sake of clarity, here the latter is referred to as ``MGS-GMRES'', as the modified Gram-Schmidt process is applied in~\cite{Saad1986} to orthogonalize successive basis vectors of~\eqref{eq:krylov}.
Later, the distinction between the two terms shall often be blurred because in the literature GMRES is almost exclusively implemented by means of a Gram-Schmidt-like process.

\subsection{MGS-GMRES}

Saad and Schultz proposed the first GMRES algorithm in~\cite{Saad1986}, which uses the Arnoldi process~\cite{Arnoldi1951} to generate an orthonormal basis~$\{v_1,\dots,v_n\}$ of the Krylov subspace~\eqref{eq:krylov}.
The key relation used in the Arnoldi process is
\begin{equation}\label{eq:cgs1}
h_{j+1,j}v_{j+1} = Av_j - \sum_{i=1}^j h_{i,j}v_i,
\end{equation}
where $v_1=r_0/\norm{r_0}$ and $h_{i,j}$ are selected such that~$v_{j+1}$ is normalized and orthogonal to all the previous basis vectors, that is,
\begin{align}
& h_{i,j} = (Av_j,v_i),\quad i=1,\dots,j, \label{eq:cgs2} \\
& h_{j+1,j} = \lVert Av_j - \sum_{i=1}^j h_{i,j}v_i \rVert, \label{eq:cgs3}
\end{align}
where $(\cdot,\cdot)$ denotes the Euclidean inner product.
The above rules correspond directly to the classical Gram-Schmidt (CGS) orthogonalization process, leading to the so-called CGS-Arnoldi algorithm.
In practice the computation of~$v_{j+1}$ by \eqref{eq:cgs1}--\eqref{eq:cgs3} undergoes a severe cancellation.
Special precautions, such as the reorthogonalization~\cite{Daniel1976,Hernandez2007} and modified Gram-Schmidt (MGS) processes, have been recommended to ensure numerical stability.
The former is costly in computation and generally used in the context of parallel algorithms~\cite{Hoemmen2010,Swirydowicz2021}, eigenvalue problems~\cite{Saad2011}, and mixed precision~\cite{Lindquist2020}, while MGS has the same operation count as CGS, but with better numerical properties; see~\cite{Leon2013} for a thorough presentation of Gram-Schmidt orthogonalization.
Arnoldi with MGS orthogonalization (MGS-Arnoldi)~\cite{Saad1980} can be viewed as the standard version of Arnoldi iteration.
\begin{algorithm}
\caption{MGS-Arnoldi\label{alg:arnoldi:mgs}}
\begin{algorithmic}[1]
  \State $v_1=r_0/\norm{r_0}$
  \For {$j=1,2,\dots,n$}
    \State $w_j=Av_j$
    \For {$i=1,2,\dots,j$}
      \State $h_{i,j}=(w_j,v_i),\quad w_j=w_j-h_{i,j}v_i$
    \EndFor
    \State $h_{j+1,j}=\norm{w_j}$,\quad $v_{j+1}=w_j/h_{j+1,j}$
  \EndFor
\end{algorithmic}
\end{algorithm}
In Algorithm~\ref{alg:arnoldi:mgs}, we can see that lines~3--6 are mathematically equivalent to~\eqref{eq:cgs1}--\eqref{eq:cgs2}.
Let $\hat{P}_j$ denote the MGS projector.
Then, according to Algorithm~\ref{alg:arnoldi:mgs}, one finds~$\hat{P}_j=(I-v_jv_j^\intercal)\dots(I-v_1v_1^\intercal)$.
Recalling the grade~$d$ of~$r_0$ defined in the preceding section, we observe that $h_{j+1,j}=0$ (line~7) if and only if~$j=d$, that is, the maximal dimension of the Krylov subspace is attained.
In what follows we shall assume~$n < d$.

Let us write
\[
V_n=[v_1,\dots,v_n]\in\mathbb{R}^{N\times n},\quad H_n=[h_{i,j}]\in\mathbb{R}^{n\times n},\quad \overbar{H}_n=
\left[\begin{array}{c}
H_n \\
h_{n+1,n}e_n^\intercal
\end{array}\right]\in\mathbb{R}^{(n+1)\times n},
\]
where $e_i$ denotes the $i$th column of the identity matrix of appropriate order and $h_{i,j}=0$ for~$i>j+1$.
Therefore, $\overbar{H}_n$ and $H_n$ are upper Hessenberg matrices.
It follows that
\begin{subequations}\label{eq:arnoldi}
\begin{align}
AV_n &= V_{n+1}\overbar{H}_n \label{eq:arnoldi:1} \\
&= V_nH_n + h_{n+1,n}v_{n+1}e_n^\intercal, \\
V_n^\mathsf{T} AV_n &= H_n,
\end{align}
\end{subequations}
which can be interpreted as an orthogonal reduction of~$A$ to upper Hessenberg form.
Let $y$ denote the coordinate vector of the correction~$x_n-x_0$ in the basis~$V_n$.
Then it follows from the minimal residual process~\eqref{eq:proj}--\eqref{eq:mr} that
\begin{align}
y &= \argmin_{\hat{y}}\norm{r_0-AV_n \hat{y}} = \argmin_{\hat{y}}\norm{\beta e_1-\overbar{H}_n \hat{y}}, \label{eq:y} \\
x_n &= x_0 + V_n y, \nonumber
\end{align}
where $\beta=\norm{r_0}$.
Saad and Schultz~\cite{Saad1986} suggested using Givens rotations that reduce the upper Hessenberg matrix~$\overbar{H}_n$ recursively to the upper triangular matrix.
\begin{algorithm}
\caption{MGS-GMRES\label{alg:gmres:mgs}}
\begin{algorithmic}[1]
  \State $\beta=\norm{r_0},\quad v_1=r_0/\beta,\quad g=\beta e_1$
  \For {$j=1,2,\dots,n$}
    \State $w_j=Av_j$
    \For {$i=1,2,\dots,j$}
      \State $h_{i,j}=(w_j,v_i),\quad w_j=w_j-h_{i,j}v_i$
    \EndFor
    \State $h_{j+1,j}=\norm{w_j},\quad v_{j+1}=w_j/h_{j+1,j}$
    \For {$i=1,2,\dots,j-1$}
      \State $\gamma=c_ih_{i,j}+s_ih_{i+1,j},\quad h_{i+1,j}=-s_ih_{i,j}+c_ih_{i+1,j},\quad h_{i,j}=\gamma$
    \EndFor
    \State $\delta=\sqrt{h_{j,j}^2+h_{j+1,j}^2},\quad s_j=h_{j+1,j}/\delta,\quad c_j=h_{j,j}/\delta$
    \State $h_{j,j}=c_jh_{j,j}+s_jh_{j+1,j},\quad h_{j+1,j}=0$
    \State $g_{j+1}=-s_jg_j,\quad g_j=c_jg_j,\quad \rho=\abs{g_{j+1}}$,\quad \textbf{if} $\rho$ small enough, \textbf{then} $n=j$, go to 15
  \EndFor
  \State Solve upper triangular linear system~$H_ny=g_{1:n}$
  \State $x_n=x_0+V_ny$
\end{algorithmic}
\end{algorithm}
We show MGS-GMRES as Algorithm~\ref{alg:gmres:mgs}, where $\rho$ denotes the residual norm, $g_i$ denotes the $i$th entry of vector~$g$ and $g_{i:j}$ refers to the subvector with elements indexed by~$i$ through~$j$.
Lines~3--7 are the same as in Algorithm~\ref{alg:arnoldi:mgs}, followed by Givens rotations
\[
\left(\begin{array}{c}h_{i,j} \\ h_{i+1,j}\end{array}\right) \leftarrow \left(\begin{array}{cc}c_i & s_i \\ -s_i & c_i\end{array}\right)\left(\begin{array}{c}h_{i,j} \\ h_{i+1,j}\end{array}\right)
\]
to maintain a QR factorization of the upper Hessenberg matrix~$\overbar{H}_j$, which involves all the previous rotations (lines~8--10) and a new rotation (lines~11--13) to annihilate the subdiagonal entry.
In particular, $s_j$ and $c_j$ correspond to the sine and cosine of the $j$th rotation angle such that $h_{j+1,j}=0$; we refer the reader to~\cite{Bindel2002} for a more robust implementation of Givens rotations.
An important by-product is that the residual norm~$\abs{g_{j+1}}$ can be readily obtained at the end of each step (line~13).
Note that the CGS-based GMRES (CGS-GMRES) algorithm can be easily derived by moving the vector update operations in line~5 out of the for-loop.
The CGS projector can be written as~$\tilde{P}_j=I-V_jV_j^\mathsf{T}$.

In~\cite{Saad1986}, it was proved that MGS-GMRES never breaks down.
Also, in a remarkable work by Paige et al.~\cite{Paige2006}, it was proved that MGS-GMRES is backward stable.
By modifying the constraint~\eqref{eq:proj2} as~$r_n \perp \mathcal{K}_n$, we can get the so-called FOM method~\cite{Saad1981}.
A detailed account of relations between GMRES and FOM can be found in~\cite{Brown1991,Cullum1996}; see also~\cite{Saad2003}.
FOM has received far less attention partly because the corresponding algorithms can break down.
However, it was shown in~\cite{Brown1991} that for the cases FOM performs poorly, GMRES will not offer much help.
In this situation, the latter may still be preferred due to its minimal residual property.
GMRES has also been compared with CMRH and Lanczos-based methods; see~\cite{Sadok1999,Sadok2012,Hochbruck1998}.

Algorithm~\ref{alg:gmres:mgs} uses long recurrences to orthogonalize~$Av_j$ against all the previous basis vectors, the cost of which grows quadratically with the number of steps.
Saad and Schultz~\cite{Saad1986} considered restarting and truncation as remedies.
The truncation strategy keeps only a small number of previous basis vectors in the orthogonalization phase, and computes the approximate solution in a progressive manner; see~\cite{Saad1996}.
In the restarted version, a parameter is introduced to limit the size of the Krylov subspaces.
The resulting scheme outlined in Algorithm~\ref{alg:gmres(m)} is called restarted GMRES, denoted by GMRES($m$).
Practically we use the restarted scheme on the top level, and a specific GMRES algorithm at the lower level; see~\cite{Fraysse2005} for an efficient implementation with variable orthogonalization schemes.
\begin{algorithm}
\caption{GMRES($m$)\label{alg:gmres(m)}}
\begin{algorithmic}[1]
  \For {$k=0,1,\dots$}
    \State Compute an $m$-dimensional basis~$Z_m$ and the coordinate vector~$y$ using a GMRES algorithm
    \State $x_{k+1}=x_k+Z_my$
    \State Compute the residual norm~$\rho$,\quad \textbf{if} $\rho$ small enough, \textbf{then} stop
  \EndFor
\end{algorithmic}
\end{algorithm}
Unfortunately, the desired properties such as finite termination and residual minimization are lost.
Both strategies may result in poor convergence or even stagnation; see,~e.g.,~\cite{Faber1996,Sosonkina1998,deSturler1999,Baker2005,Giraud2010,Soodhalter2020}.
Moreover, using a large~$m$ may not necessarily lead to fast convergence (see~\cite{Eiermann2000,Embree2003}).
We shall discuss other potential acceleration strategies in Section~\ref{sec:acc}.

\subsection{HH-GMRES}

The minimal residual Krylov subspace method can be completely described by~\eqref{eq:proj} from a mathematical point of view.
From a computational point of view, however, mathematically equivalent algorithms may have quite different numerical behavior.
Since at that time the backward stability of MGS-GMRES~\cite{Paige2006} was still unknown, motivated by numerical stability concerns, Walker~\cite{Walker1988} in 1988 proposed an algorithm that uses Householder reflections to orthogonalize the basis vectors; see also~\cite{Walker1989}.
\begin{algorithm}
\caption{HH-GMRES\label{alg:gmres:hh}}
\begin{algorithmic}[1]
  \State $w_1=r_0+\sign(r_{0,1})\norm{r_0}e_1$,\quad $P_1=I-2w_1w_1^\intercal/(w_1^\intercal w_1)$,\quad $\beta=-\sign(r_{0,1})\norm{r_0}$,\quad $v_1=P_1 e_1$
  \For {$j=1,2,\dots,n$}
    \State $u=P_jP_{j-1}\dots P_1 Av_j$,\quad $w_{j+1}=[0,\dots,0,u_{j+1},\dots,u_N]^\mathsf{T}$
    \State $w_{j+1}=w_{j+1}+\sign(u_{j+1})\norm{w_{j+1}}e_{j+1}$,\quad $P_{j+1}=I-2w_{j+1}w_{j+1}^\intercal/(w_{j+1}^\intercal w_{j+1})$
    \State $h_j=P_{j+1}u$,\quad $v_{j+1}=P_1P_2\dots P_{j+1}e_{j+1}$
  \EndFor
  \State Define~$\overbar{H}_n$ as the first~$n+1$ rows of the matrix~$[h_1,\dots,h_n]$
  \State Compute~$y$ such that $\norm{\beta e_1-\overbar{H}_n y}$ is minimized (see lines~8--15 of Algorithm~\ref{alg:gmres:mgs})
  \State $x_n=x_0+V_ny$
\end{algorithmic}
\end{algorithm}
Algorithm~\ref{alg:gmres:hh} depicts the Householder-based GMRES (HH-GMRES) algorithm.
Here, $\sign(\alpha)$ returns~$1$ or~$-1$ depending on the sign of~$\alpha$ and $r_{i,j}$ denotes the~$j$th element of~$r_i$.
The signs (lines~1 and~4) are chosen to reduce the risk of subtractive cancellation.
Line~4 corresponds directly to Householder reflections; the details can be found in~\cite{Golub2013,Saad2003}.
HH-GMRES gradually generates a QR factorization
\[
[r_0,Av_1,\dots,Av_n] = P_1P_2\dots P_{n+1}[\beta e_1,h_1,\dots,h_n].
\]
By discarding the zero rows of the upper triangular matrix and the first column of both sides, we can obtain the Arnoldi relation~\eqref{eq:arnoldi}.
The Householder vector~$w_{j+1}$ contains $j$ zeros in the first rows (line~3) such that the Householder transformation can leave the previous vectors~$h_1,\dots,h_{j-1}$ unchanged.
For the solution of line~8 and the estimate of the residual norm at each iteration step, the same techniques used in Algorithm~\ref{alg:gmres:mgs} can be applied.
Note that we only need to store the~$w_i$'s instead of the~$P_i$'s.
We also prefer not to explicitly store the~$v_i$'s but use the Qin Jiushao's scheme (Horner's scheme) to evaluate line~9 in a recursive manner; see,~e.g.,~\cite{Saad2003}.

Similar to MGS-GMRES, the cost of HH-GMRES increases dramatically as the number of steps increases.
Therefore, the restarted scheme in Algorithm~\ref{alg:gmres(m)} is commonly employed as an outer loop.
HH-GMRES was proved earlier to be backward stable (see~\cite{Drkosova1995}), but can be roughly two times more expensive than MGS-GMRES.

\subsection{Simpler GMRES}

In~1994, Walker and Zhou~\cite{Walker1994} proposed another variant of the GMRES method, called simpler GMRES (SGMRES), which does not require the QR factorization of an upper Hessenberg matrix. 
Although rarely used in practice due to stability issues, this algorithm can enhance our understanding of minimal residual iterations.

We present a generalized version as given in~\cite{Jiranek2008} that contains the essential features of SGMRES.
Consider a nonorthogonal basis~$Z_n=\{z_1,\dots,z_n\}$ of~$\mathcal{K}_n$ and an orthonormal basis~$V_n=\{v_1\dots,v_n\}$ of~$A\mathcal{K}_n$.
We construct a matrix relation
\begin{equation}\label{eq:sgmres}
AZ_n=V_nT_n
\end{equation}
such that $T_n=[t_{i,j}]\in\mathbb{R}^{n\times n}$ is an upper triangular matrix.
According to the constraint~\eqref{eq:proj2}, one finds $r_n=(I-V_nV_n^\mathsf{T})r_0$, yielding the following update:
\begin{equation}\label{eq:sgmres:r}
\alpha_n=(r_{n-1},v_n),\quad r_n=r_{n-1}-\alpha_nv_n.
\end{equation}
In addition, since~$x_n=x_0+Z_ny$ with some n-dimensional vector~$y$, it follows that~$r_n=r_0-V_nT_ny$.
Combining~\eqref{eq:proj2} and~\eqref{eq:sgmres:r}, one obtains that
\[
T_ny=V_n^\mathsf{T} r_0=[\alpha_1,\dots,\alpha_n]^\mathsf{T}.
\]
\begin{algorithm}
\caption{Generalized Simpler GMRES\label{alg:gsgmres}}
\begin{algorithmic}[1]
  \For {$j=1,2,\dots,n$}
    \State Choose~$z_j$ and compute~$w_j=Az_j$
    \State Orthonormalize~$w_j$ against~$v_1,\dots,v_{j-1}$ to obtain $t_{1,j},\dots,t_{j,j}$ and~$v_j$
    \State $\alpha_j=(r_{j-1},v_j),\quad r_j=r_{j-1}-\alpha_jv_j$
    \State Compute the residual norm~$\rho$,\quad \textbf{if} $\rho$ small enough, \textbf{then} $n=j$, go to 7
  \EndFor
  \State Solve $T_ny=[\alpha_1,\dots,\alpha_n]^\mathsf{T}$
  \State $x_n=x_0+Z_ny$
\end{algorithmic}
\end{algorithm}
The generalized simpler GMRES framework (see~\cite{Jiranek2008}) is listed in Algorithm~\ref{alg:gsgmres}.
Choosing~$Z_n=[r_0/\norm{r_0},V_{n-1}]$ results in SGMRES, which was formalized by virtue of both MGS- and Householder-based processes in~\cite{Walker1994}.
Note that in the Householder variant there is no need to store the~$v_i$'s and thus line~4 should be modified accordingly; see~\cite{Walker1994} for details.
The choice~$Z_n=[r_0/\norm{r_0},\dots,r_{n-1}/\norm{r_{n-1}}]$ was introduced in~\cite{Jiranek2008}, leading to a slightly different algorithm called residual-based simpler GMRES (RB-SGMRES).
Then, an adaptive variant combining the two choices was proposed in~\cite{Jiranek2010}, following the rule
\begin{equation}\label{eq:asgmres}
z_n =
\begin{cases}
\frac{r_{n-1}}{\norm{r_{n-1}}},& \norm{r_{n-1}}\le\omega\norm{r_{n-2}}, \\
v_{n-1}, & \text{otherwise},
\end{cases}
\end{equation}
with~$\omega\in[0,1]$.
Notice that $\omega=0$ and $\omega=1$ correspond to SGMRES and RB-SGMRES, respectively.

Theoretical aspects are collected in~\cite{Walker1994,Liesen2000,Liesen2002,Chen2004,Jiranek2008,Jiranek2010}, which provides bounds for the condition number of~$Z_n$ and allows for a better understanding of the unstable behavior of Algorithm~\ref{alg:gsgmres}; see also~\cite{Meurant2020}.
It's interesting to observe that, along with the results in~\cite{Greenbaum1994}, SGMRES can be employed to describe GMRES convergence; see~\cite{Liesen2000}.

\subsection{Equivalent formulations}\label{sec:ef}

We briefly review some iterative schemes that are in some sense mathematically equivalent to GMRES.
Paige and Saunders~\cite{Paige1975} in 1975 proposed the MINRES algorithm satisfying~\eqref{eq:mr} for solving symmetric indefinite systems, which motivated the development of MGS-GMRES.
In 1982, Elman~\cite{Elman1982,Eisenstat1983} proposed the GCR algorithm which amounted to generalizing the conjugate residual algorithm (see,~e.g.,~\cite{Saad2003}) to the nonsymmetric case.
\begin{algorithm}
\caption{GCR\label{alg:gcr}}
\begin{algorithmic}[1]
  \State $q_0=r_0$
  \For {$j=0,1,\dots,n$}
    \State $\alpha_j=(r_j,Aq_j)/(Aq_j,Aq_j)$
    \State $x_{j+1}=x_j+\alpha_jq_j,\quad r_{j+1}=r_j-\alpha_jAq_j$
    \State Compute the residual norm~$\rho$,\quad \textbf{if} $\rho$ small enough, \textbf{then} stop
    \State $\beta_{i,j}=-(Ar_{j+1},Aq_i)/(Aq_i,Aq_i),\quad q_{j+1}=r_{j+1}+\sum_{i=0}^j\beta_{i,j}q_i$
  \EndFor
\end{algorithmic}
\end{algorithm}
In Algorithm~\ref{alg:gcr}, the term~$Aq_j$ can be updated by means of the formula in line~6 and the value of~$Ar_j$ without an additional matrix-vector multiplication.
GCR generates a set of $A$-orthogonal basis vectors~$q_j$ spanning the Krylov subspace~\eqref{eq:krylov}.
These vectors can be alternatively updated by
\begin{equation}\label{eq:orthodir}
\beta_{i,j}=-(A^2q_j,Aq_i)/(Aq_i,Aq_i),\quad q_{j+1}=Aq_j+\sum_{i=0}^j\beta_{i,j}q_i,
\end{equation}
leading to the so-called ORTHODIR algorithm, which was proposed by Jea and Young~\cite{Young1980,Jea1983}.
A truncated version of GCR was given earlier in~\cite{Vinsome1976}, and both GCR and ORTHODIR can be formulated in a cyclic manner like Algorithm~\ref{alg:gmres(m)} to reduce the associated work and storage.
Nonetheless, GCR can fail if the symmetric part of~$A$ is indefinite and the implementation of~\eqref{eq:orthodir} can cause stability problems, while it is known that MGS-GMRES is a cheaper and more robust alternative; see,~e.g.,~\cite{Saad1986}.
Similar reasoning holds for another early variant proposed by Axelsson~\cite{Axelsson1980}.
It was mentioned in~\cite{Walker1994} that there is a close relationship between SGMRES and ORTHODIR; see also~\cite{Jiranek2008}, in which it was also shown that RB-SGMRES is closely related to GCR.

An early algorithm that realized~\eqref{eq:mr} was already described by Khabaza~\cite{Khabaza1963} in 1963, which used $r_0,Ar_0,\dots,A^{n-1}r_0$ as the basis vectors and thus suffered from numerical instabilities.
In~\cite{Eirola1989}, Eirola and Nevanlinna formalized a quasi-Newton method, later called EN method, that approximates~$A^{-1}$ by rank-one updates.
It was shown in~\cite{Eirola1989,Vuik1992} that a variant of EN is mathematically equivalent to GMRES; see also~\cite{vanderVorst1994} for additional comments.
More recent examples of equivalence relations involving quasi-Newton methods can be found in~\cite{Haelterman2010,Haelterman2015}.
In 1988, on the other hand, Sidi gave a proof in~\cite{Sidi1988} that an acceleration technique called reduced rank extrapolation is mathematically equivalent to GMRES; see~\cite{Smith1987} for a survey on vector extrapolation methods.
Later, it was established by Walker and Ni~\cite{Walker2011} that another technique called Anderson acceleration is equivalent to GMRES in some sense; we refer the reader to~\cite{Walker2011,Brezinski2018} for more details on Anderson acceleration and its relation to GMRES.

\subsection{Convergence behavior}

An important observation for studying GMRES convergence is that any vector in the Krylov subspace~\eqref{eq:krylov} can be written in polynomial form, which implies
\[
r_n = p_n(A)r_0,\quad p_n\in\mathcal{P}_n,
\]
where $\mathcal{P}_n$ is the set of all polynomials~$p$ of degree at most~$n$ such that $p(0)=1$.
Then, \eqref{eq:mr} can be expressed as
\begin{equation}\label{eq:mr:poly}
\norm{r_n} = \norm{p_n(A)r_0} = \min_{p\in\mathcal{P}_n}\norm{p(A)r_0}.
\end{equation}
Let $\lambda_i(\cdot)$ denote the eigenvalues of a matrix and $\lambda_i=\lambda_i(A)$.
Assume that $A$ is diagonalizable so that~$A=X\Lambda X^{-1}$ with $\Lambda=\diag(\lambda_1,\dots,\lambda_N)$.
The residual norm of GMRES satisfies
\begin{equation}\label{eq:bd:elman1982}
\frac{\norm{r_n}}{\norm{r_0}} \le \kappa(X)\min_{p\in\mathcal{P}_n}\max_i\abs{p(\lambda_i)},
\end{equation}
where $\kappa(X)=\norm{X}\norm{X^{-1}}$ denotes the condition number of~$X$.
Note that $A$ being normal implies~$\kappa(X)=1$, in which case the bound in~\eqref{eq:bd:elman1982} is sharp (see~\cite{Greenbaum1994c,Joubert1994}).
This bound was presented in 1982 by Elman~\cite{Elman1982} in the context of the so-called GCR algorithm; see also~\cite{Eisenstat1983,Saad1986}.
It was also proved in~\cite{Elman1982} that if $M=(A+A^\mathsf{T})/2$ is positive definite, then we have
\[
\frac{\norm{r_n}}{\norm{r_0}} \le \left(1-\frac{{\lambda_{\min}(M)}^2}{\lambda_{\max}(A^\mathsf{T} A)}\right)^{\frac{n}{2}},
\]
where $\lambda_{\min}(\cdot)$ and $\lambda_{\max}(\cdot)$ are respectively the smallest and largest eigenvalues of a matrix in the absolute sense.
This bound was improved two decades later in~\cite{Beckermann2005b,Beckermann2005}.
Now, it follows from~\eqref{eq:mr:poly} that
\begin{equation}\label{eq:bd:worst}
\frac{\norm{r_n}}{\norm{r_0}} \le \psi_n(A) \le \varphi_n(A),\quad \psi_n(A)=\max_{\norm{v}=1}\min_{p\in\mathcal{P}_n}\norm{p(A)v},\quad \varphi_n(A)=\min_{p\in\mathcal{P}_n}\norm{p(A)}.
\end{equation}
In general, $\psi_n(A)$ is called worst-case GMRES value, as its purpose is to study the worst-case behavior of the GMRES method.
It is known that $\psi_n(A)$ is attainable by the residual norm, that is, for each~$n$ there exists an initial vector~$v$ such that $\norm{r_n}=\psi_n(A)$ (see,~e.g.,~\cite{Tichy2007}).
$\varphi_n(A)$ is called ideal GMRES value, which was introduced in~\cite{Greenbaum1994b} to totally exclude the influence of the initial vector.
It was proved in~\cite{Greenbaum1994b,Liesen2009} that $\varphi_n(A)$ has a unique minimizer, while as discussed in~\cite{Faber2013,Liesen2004c}, the polynomial and unit-norm vector associated with~$\psi_n(A)$ may not be uniquely determined.
Moreover, if $A$ is normal, then~$\psi_n(A)=\varphi_n(A)$ holds (see~\cite{Greenbaum1994c,Joubert1994}); otherwise, there exist examples for which~$\psi_n(A)<\varphi_n(A)$ (see~\cite{Faber1996,Toh1997,Faber2013}).
We refer the reader to~\cite{Liesen2004,Liesen2004c,Tichy2007,Faber2013} and their references for further details on~\eqref{eq:bd:worst}.

The field of values plays a beneficial role in convergence analysis, which is defined as
\[
\mathcal{F}(A) = \{(Av,v):\norm{v}=1,v\in\mathbb{C}^N\}.
\]
Assume that the origin is outside~$\mathcal{F}(A)$ and let
\[
\mu_F(A) = \min_{\lambda\in\mathcal{F}(A)}\abs{\lambda}.
\]
Then, the following bound holds:
\begin{equation}\label{eq:bd:field}
\frac{\norm{r_n}}{\norm{r_0}} \le (1-\mu_F(A)\mu_F(A^{-1}))^{\frac{n}{2}};
\end{equation}
see~\cite{Starke1997,Eiermann2001} for proofs based on the worst-case GMRES approximation; see also~\cite{Liesen2012} where the bound in~\eqref{eq:bd:field} was proved to hold for the ideal GMRES approximation.
Besides, an early account of the field of the values applied to convergence analysis can be found in~\cite{Eiermann1993}.
We refer the reader to the most recent work~\cite{Benzi2021} and the references therein for the use of the field of values for convergence analysis.
The polynomial numerical hull introduced in~\cite{Nevanlinna1993} can be seen as a generalization of the field of values, which is defined as
\[
\mathcal{H}_n(A) = \{\lambda\in\mathbb{C}:\norm{p(A)}\ge\abs{p(\lambda)},\deg(p)\le n\},
\]
from which a lower bound of the ideal GMRES approximation can be easily derived.
However, the determination of~$\mathcal{H}_n(A)$ is very difficult in its full generality.
We refer to~\cite{Greenbaum2002,Faber2003,Greenbaum2004,Tichy2007} for related works.
Another group of results is obtained via the $\epsilon$-pseudospectrum~\cite{Trefethen1990}.
Let $I$ denote the identity matrix of appropriate size.
For a real $\epsilon>0$, the $\epsilon$-pseudospectrum of~$A$ is the set
\[
\Lambda_\epsilon = \{\lambda\in\mathbb{C}:\norm{\lambda I - A}^{-1} \ge \epsilon^{-1}\},
\]
from which an upper bound of~$\varphi_n(A)$ can be obtained; see,~e.g.,~\cite{Nachtigal1992,Trefethen2005} for more details.

Linear contraction bounds for the residual norms can be somewhat misleading since a highly nonlinear convergence behavior can not be adequately described by a linear contraction.
The book~\cite{Liesen2013} written by Liesen and Strako{\v{s}} is very useful for getting more insight in the cost of computations (see also~\cite{Carson2020b}), which highlights the inherent nonlinear nature of Krylov subspace methods.
In~\cite{Embree1999}, an experimental comparison of GMRES convergence bounds based on eigenpairs, field of values and $\epsilon$-pseudospectrum for non-normal matrices can be found, but none of these can be consistently better than the others.
In order to get more descriptive results, Titley-Peloquin et al.~\cite{TPeloquin2014} in 2014 derived several GMRES bounds that involve the initial residual vector; see also~\cite{Pestana2013} for an analysis with respect to nonstandard inner products.
More recently, Sacchi and Simoncini~\cite{Sacchi2019} gave a more descriptive convergence analysis specifically for the case of localized ill-conditioning.

It was confirmed in a series of papers by Arioli, Greenbaum, Pt{\'a}k and Strako{\v{s}}~\cite{Greenbaum1994,Greenbaum1996,Arioli1998} from 1994 to 1998 that convergence of GMRES for non-normal matrices could not be determined solely by the distribution of eigenvalues.
In~\cite{Greenbaum1994}, Greenbaum and Strako{\v{s}} studied the matrices~$B$ for which the spaces $A\mathcal{K}_n(A,r_0)$ and $B\mathcal{K}_n(B,r_0)$ are the same.
As a result, the convergence history generated for the $(B,r_0)$ pair by GMRES is the same as that generated for the $(A,r_0)$ pair.
The matrices~$B$ that satisfies this property are called GMRES-equivalent matrices (see,~e.g.,~\cite{Liesen2000,DTebbens2014b}).
Among other results, Greenbaum and Strako{\v{s}} concisely state in~\cite{Greenbaum1994} on GMRES convergence that
\begin{quote}
Any behavior that can be seen with the method can be seen with the method applied to a matrix having any nonzero eigenvalues.
\end{quote}
This has been extended by Greenbaum et al.~\cite{Greenbaum1996} along a different line, who conclude that
\begin{quote}
Any nonincreasing convergence curve can be obtained with GMRES applied to a matrix having any desired eigenvalues.
\end{quote}
Finally, a complete parametrization of all~$(A,b)$ pairs with prescribed eigenvalues for which GMRES generates prescribed residual norms has been given in~\cite{Arioli1998}; see also~\cite{Meurant2012} for further discussion.
GMRES-equivalent matrices have been employed by some subsequent studies~\cite{Liesen2000,DTebbens2014b}, while investigations of ``any behavior is possible'', like that in~\cite{Greenbaum1994,Greenbaum1996,Arioli1998}, have been continued in the past ten years, concerning Ritz values~\cite{DTebbens2012,DTebbens2014}, harmonic Ritz values~\cite{Du2017}, restarted GMRES~\cite{Vecharynski2011,DTebbens2020}, and block GMRES~\cite{Kubinova2020}.
Note that Ritz and harmonic Ritz values are respectively the roots of FOM and GMRES residual polynomials; see~\cite{Goossens1999,Meurant2017}.
We also mention~\cite{Meurant2012} for a generalization to other Krylov subspace methods.

Eigenvalue clusters associated with a parameterized model were employed in~\cite{Campbell1996} to explain GMRES convergence.
Exact but complicated expressions for residual norms were derived in~\cite{Ipsen2000,Liesen2004b,Sadok2005}; see also~\cite{Meurant2020}.
Inspired by the work~\cite{vanderSluis1986} on the conjugate gradient method, van der Vorst and Vuik~\cite{vanderVorst1993} in 1993 investigated the superlinear convergence behavior of GMRES by means of Ritz values.
Simoncini and Szyld~\cite{Simoncini2005} addressed the same problem using spectral projectors.
Moret~\cite{Moret1997} discussed superlinear convergence of GMRES in a complex separable Hilbert space; see also~\cite{Blechta2021} for a generalization of Moret's result.
Since any nonincreasing convergence curve is possible for GMRES~\cite{Greenbaum1996}, on the other hand, stagnation may occur even for the original GMRES method without restarting.
This was first investigated by Brown~\cite{Brown1991} in 1991, and was then discussed in~\cite{Zavorin2003,Liesen2004c,Simoncini2008,Simoncini2010,Meurant2012b,Meurant2014}.
Although most work on convergence analysis considered full GMRES, the restarted version is more practical, on which we refer the reader to,~e.g.,~\cite{Joubert1994b,Zitko2000,Zitko2008,Baker2009,Vecharynski2010,Vecharynski2011,DTebbens2020} for more details.

The above discussion assumes exact arithmetic.
In finite precision arithmetic, it was proved by Drko{\v{s}}ov{\'a} et al.~\cite{Drkosova1995} in 1995 that HH-GMRES is backward stable.
In 2006, it was proved by Paige et al.~\cite{Paige2006} that MGS-GMRES is also backward stable; previous observations can be found in~\cite{Greenbaum1997c,Paige2002}.
Due to the ill-conditioning of the matrix~$Z_n$ in~\eqref{eq:sgmres}, SGMRES has serious stability problems (see,~e.g,~\cite{Liesen2002,Chen2004}), which can be partially remedied by using RB-SGMRES~\cite{Jiranek2008} or the adaptive variant~\eqref{eq:asgmres}~\cite{Jiranek2010}; see also the experiments in~\cite{Matinfar2012}.

\section{Acceleration strategies}\label{sec:acc}

The focus in this section is on the acceleration strategies to tackle the challenge that a basic GMRES algorithm might face.
First, it is almost certain that the use of preconditioners will be beneficial for Krylov subspace methods when solving ill-conditioned problems.
There are also numerous publications concerning deflation and augmentation, some of which are highly related to preconditioning.
Then we briefly introduce weighted inner products and inexact matrix-vector products, and end up this section by discussing mixed-precision techniques.

\subsection{Preconditioning}\label{sec:prec}

Preconditioning of the linear system is used to improve the convergence behavior of Krylov subspace methods.
One needs to find a nonsingular matrix~$M$ close to~$A$ in some sense, commonly called preconditioning matrix or preconditioner, such that the new system involving~$M$ should be inexpensive to solve, that is,
\begin{align}
& M^{-1}Ax=M^{-1}b, \label{eq:pre:left} \\
& AM^{-1}u=b,\quad Mx=u. \label{eq:pre:right}
\end{align}
The cases~\eqref{eq:pre:left} and~\eqref{eq:pre:right} correspond to left preconditioning and right preconditioning, respectively.
Another form called split preconditioning can be employed when nearly symmetric systems are encountered; see,~e.g.,~\cite{Saad2003}.
For GMRES, right preconditioning is often preferred, because the residuals in this case are identical to the true residuals, while left preconditioning leads to preconditioned residuals~$z_n=M^{-1}(b-Ax_n)$.
One also needs to transform line~3 in Algorithm~\ref{alg:gmres:mgs} to~$w_j=M^{-1}Av_j$ for the case~\eqref{eq:pre:left} and~$w_j=AM^{-1}v_j$ for~\eqref{eq:pre:right}.
High-quality preconditioners are often constructed by incorporating problem-specific information; see,~e.g.,~\cite{Nigro1998,Erlangga2008b,Erlangga2008c,Gander2015b,Pearson2020}.
On the other hand, general-purpose preconditioners have also been broadly exploited, such as incomplete factorization, approximate inverse, and algebraic multigrid; see~\cite{Benzi2002,Wathen2015} and references therein.

The flexible MGS-GMRES variant proposed by Saad~\cite{Saad1993} in 1993 allows the preconditioner to vary at each step, which shall be denoted by FGMRES as widely used in the literature.
\begin{algorithm}
\caption{FGMRES\label{alg:fgmres}}
\begin{algorithmic}[1]
  \State $\beta=\norm{r_0},\quad v_1=r_0/\beta$
  \For {$j=1,2,\dots,n$}
    \State $z_j=M_j^{-1}v_j,\quad w_j=Az_j$
    \For {$i=1,2,\dots,j$}
      \State $h_{i,j}=(w_j,v_i),\quad w_j=w_j-h_{i,j}v_i$
    \EndFor
    \State $h_{j+1,j}=\norm{w_j},\quad v_{j+1}=w_j/h_{j+1,j}$
  \EndFor
  \State Compute~$y$ such that $\norm{\beta e_1-\overbar{H}_n y}$ is minimized
  \State $x_n=x_0+Z_ny$
\end{algorithmic}
\end{algorithm}
As seen in Algorithm~\ref{alg:fgmres}, $z_j=M_j^{-1}v_j$ should be explicitly stored in~$Z_n=[z_1,\dots,z_n]$, and then reused in line~10.
Therefore, the Arnoldi-like relation~$AZ_n=V_{n+1}\overbar{H}_n$ holds.
Unlike standard right preconditioning, which constructs an orthonormal basis of the Krylov subspace
\[
\vspan\{r_0,AM^{-1}r_0,\dots,(AM^{-1})^{n-1}r_0\},
\]
the subspace of solution estimates in FGMRES is no longer a standard Krylov subspace.
In addition, $H_n$ may be singular even if $A$ is nonsingular, and thus the assumption of nonsingularity of~$H_n$ must be made when deducing exact solution from~$h_{j+1,j}=0$; see~\cite{Saad1993,Saad2003} for more discussion.
A very detailed analysis of FGMRES, in comparison with right preconditioned MGS-GMRES, can be found in~\cite{Arioli2007}.
More recently, Greif et al.~\cite{Greif2017} developed the multi-preconditioned GMRES scheme in 2017, which shares some similarities with the block version of FGMRES~\cite{Calandra2012,Calandra2013}, but targets~\eqref{eq:ls} instead of systems with multiple right-hand sides.
In practice, a restarted procedure should be used in replacement of the simple version in Algorithm~\ref{alg:fgmres}; see~\cite{Fraysse2008} for details on the implementation of FGMRES.

Another interesting but less used variable preconditioning technique proposed by van der Vorst and Vuik~\cite{vanderVorst1994} in 1994 consists of the GCR algorithm as the outer iterations and a GMRES algorithm as the inner iterations.
The inner algorithm constructs a preconditioner for the outer procedure.
This new scheme is called GMRESR, which was motivated by the EN algorithm~\cite{Eirola1989}.
Note that in GMRESR the residuals are preconditioned, while in FGMRES, as shown in Algorithm~\ref{alg:fgmres}, the search directions are preconditioned.
We refer to~\cite{Vuik1995} for a comparison of FGMRES and GMRESR.
Extensions of GMRESR have been developed in~\cite{deSturler1996,deSturler1999} by requiring that the inner basis maintains orthogonality to the outer basis, possibly combining the truncation strategy to reduce storage costs.

The GMRES polynomial itself, defined as the polynomial in~\eqref{eq:mr:poly}, can be used to derive preconditioners for Krylov subspace algorithms.
Liu et al.~\cite{Liu2015} in 2015 developed a polynomial preconditioned GMRES variant using the GMRES polynomial.
An improved approach was proposed in~\cite{Loe2022} by Loe and Morgan.
The idea can be summarized as follows:
\begin{enumerate}
\item Run $m$ steps of MGS-GMRES;
\item Compute the harmonic Ritz values (see~\eqref{eq:harmonic}) which are the roots of the GMRES polynomial;
\item Order the roots and construct the polynomial;
\item Apply a GMRES algorithm to the polynomial preconditioned system.
\end{enumerate}
We refer the reader to~\cite{Liu2015,Loe2020,Loe2022} for more details and experiments.
Polynomial preconditioners can be combined with other preconditioners and were highlighted as effective tools for parallel computing.
Early contributions of hybrid or polynomial preconditioned GMRES algorithms can be found in~\cite{Nachtigal1992b,Starke1993,Joubert1994,vanGijzen1995}.
The stability problems arising in high degree preconditioners were addressed in~\cite{Embree2021,Loe2022,Ye2021}.

In order to overcome scaling difficulties for large number of processor cores, two general iterative schemes were introduced in 2014 by McInnes et al.~\cite{McInnes2014}.
As an example, they presented a two-level hierarchical FGMRES algorithm.
The outer level applies FGMRES to the global problem, while the inner level executes inexact GMRES within subgroups of cores.
Therefore, the inner solver serves as a variable preconditioner for the outer solver and entails less communication overhead.
The two-level scheme can be extended to multiple levels with multiple solvers, leading to the hierarchical Krylov scheme.
Meanwhile, a nested version was also presented in~\cite{McInnes2014}.
The inner solvers in the nested scheme apply inexact Krylov algorithms to the global system instead of local systems.

\subsection{Deflation and augmentation}

The convergence of GMRES($m$) is often hampered by eigenvalues of small magnitude and loss of information at the end of each cycle.
Deflation techniques can eliminate certain small eigenvalues or move certain small eigenvalues away from the origin.
Note that the word ``deflation'' has multiple meanings in the literature which may lead to confusion.
Augmentation amounts to considering an $m$-dimensional subspace~$\mathcal{K}_{m_1}+\vspan\{u_1,\dots,u_{m_2}\}$ with~$m=m_1+m_2$, adding information from previous cycles to the search space.

Let us write
\begin{equation}\label{eq:zm}
Z_m=[z_1,\dots,z_m]=[v_1,\dots,v_{m_1},u_1,\dots,u_{m_2}],
\end{equation}
where the first $m_1$ vectors are the orthonormal Arnoldi vectors.
In 1995, Morgan~\cite{Morgan1995} presented an augmented GMRES algorithm, later called GMRES-E in~\cite{Morgan2000}, where $u_1,\dots,u_{m_2}$ are chosen as the harmonic Ritz vectors corresponding to the harmonic Ritz values of smallest magnitude; see also~\cite{Chapman1997,Saad1997} for more analysis of this approach.
Given the relation
\begin{equation}\label{eq:gmrese}
AZ_m=V_{m+1}\overbar{H}_m,
\end{equation}
where $V_{m+1}$ forms a set of orthonormal vectors, it is known that solving the harmonic eigenvalue problem
\[
Z_m^\mathsf{T}A^\mathsf{T}Z_m y_i = \frac{1}{\theta_i}Z_m^\mathsf{T}A^\mathsf{T}AZ_m y_i
\]
is equivalent to solving
\begin{equation}\label{eq:harmonic}
(H_m+h_{m+1,m}^2H_m^{-\mathsf{T}}e_me_m^\intercal) y_i=\theta_i y_i,
\end{equation}
with $u_i=Z_m y_i$.
Morgan~\cite{Morgan2000} in 2000 proposed a mathematically equivalent but more efficient algorithm, commonly called GMRES-IR, which uses the implicit restarting technique with unwanted harmonic Ritz values as shifts to generate augmented subspaces (see~\cite{Lehoucq1996}).
A particularly important observation is that the augmented subspace can be formulated as a Krylov subspace with a carefully chosen starting vector~$\hat{r}_0$, that is,
\[
\vspan\{\hat{r}_0,A\hat{r}_0,\dots,A^{m-1}\hat{r}_0\} = \vspan\{r_0,Ar_0,\dots,A^{m_1-1}r_0,u_1,\dots,u_{m_2}\}.
\]
We mention that Le Calvez and Molina~\cite{LeCalvez1999} developed independently an implicitly restarted GMRES algorithm, which, however, resorts to a different eigenvalue problem instead of~\eqref{eq:harmonic}.
Further, Morgan~\cite{Morgan2002} revisited GMRES-IR in 2002 and found that there was room for improvement.
The so-called GMRES with deflated restarting (GMRES-DR) algorithm~\cite{Morgan2002}, based on the thick restarting approach (see~\cite{Wu2000}), can be simpler and more stable than GMRES-IR; see also~\cite{Rollin2008} for a slightly modified variant.
Besides, there are many situations where variable preconditioners are worth to be considered, and thus a flexible variant of GMRES-DR (FGMRES-DR) was suggested in~\cite{Giraud2010}.
We point here to the more recent work of Liu et al.~\cite{Liu2015} who applied polynomial preconditioning to GMRES-DR.
In addition, the effectiveness of GMRES-DR was studied by Morgan et al.~\cite{Morgan2016}.
The main idea is that the approximate eigenvectors generated by GMRES-DR can be seen as pseudoeigenvectors.
It was shown in~\cite{Morgan2016} by examples that GMRES-DR can also work for highly nonnormal matrices.

Some deflation strategies are highly related to preconditioning.
Kharchenko and Yeremin \cite{Kharchenko1995} in 1995 suggested the use of low-rank transformations by means of left and right Ritz vectors to construct a right preconditioner.
As a result, extremal eigenvalues would be translated into a vicinity of one.
In 1996, Erhel et al.~\cite{Erhel1996} employed a different approach to estimate the invariant subspace and constructed a right preconditioner that can move the small eigenvalues to a multiple large eigenvalue.
A similar idea was later suggested by Baglama et al.~\cite{Baglama1998}.
They described two algorithms in the context of left preconditioning, which exploit implicit restarting to improve flexibility.
Moreover, an improved algorithm built on the previous work~\cite{Morgan1995,Erhel1996,Baglama1998} was introduced in~\cite{Burrage1998}. See also~\cite{Moriya2000,Baker2009} for an adaptive strategy for determining the restart frequency.

In 1996, de Sturler~\cite{deSturler1996} proposed an inner-outer algorithm with augmented basis, called GCRO, which was extended to the so-called GCROT algorithm in~\cite{deSturler1999} by incorporating truncation strategies.
A beautiful account of Morgan's work, deflation by preconditioning, and GCROT can be found in the survey~\cite{Eiermann2000}.
It was observed by Baker et al.~\cite{Baker2005} in 2005 that GMRES($m$) has alternating behavior, namely, the residual vector often alternates direction at the end of each cycle, leading to slow convergence and even stalling.
They gave some remedies to overcome this problem by combining the ideas of augmented subspaces~\cite{Morgan1995} and reserved error approximations~\cite{deSturler1996}, resulting in the LGMRES($m_1,m_2$) algorithm.
\begin{algorithm}
\caption{LGMRES($m_1,m_2$)\label{alg:lgmres}}
\begin{algorithmic}[1]
  \State $m=m_1+m_2,\quad \beta=\norm{r_0},\quad v_1=r_0/\beta$
  \For {$k=0,1,\dots$}
  \For {$j=1,2,\dots,m$}
    \If {$j\le m_1$ \textbf{or} $j-m_1>k$}
      \State $w_j=Av_j$
    \Else
      \State $w_j=Au_{k-(j-m_1-1)}$
    \EndIf
    \For {$i=1,2,\dots,j$}
      \State $h_{i,j}=(w_j,v_i),\quad w_j=w_j-h_{i,j}v_i$
    \EndFor
    \State $h_{j+1,j}=\norm{w_j},\quad v_{j+1}=w_j/h_{j+1,j}$
  \EndFor
  \State Compute~$y$ such that $\norm{\beta e_1-\overbar{H}_m y}$ is minimized
      \If {$k-m_2<0$}
      \State $Z_m=[v_1,\dots,v_{m_1+m_2-k},u_k,\dots,u_1]$
    \Else
      \State $Z_m=[v_1,\dots,v_{m_1},u_k,\dots,u_{k-m_2+1}]$
    \EndIf
  \State $u_{k+1}=Z_my,\quad x_{k+1}=x_k+u_{k+1},\quad r_{k+1}=b-Ax_{k+1}$
  \State $\beta=\norm{r_{k+1}}$,\quad \textbf{if} $\beta$ small enough, \textbf{then} stop; \textbf{else} $v_1=r_{k+1}/\beta$
  \EndFor
\end{algorithmic}
\end{algorithm}
Algorithm~\ref{alg:lgmres} satisfies~\eqref{eq:gmrese} and generalizes GMRES-E to a GCRO-like vector selection strategy.
Notice that~$u_{k+1}$ can be considered as an approximation to the exact error, which is the key idea involved in both GCRO and LGMRES; see~\cite{deSturler1996,Baker2005}.
We refer the reader to~\cite{Baker2005} for further details and numerical experiments; see also~\cite{Imakura2016} for a similar approach.
On the other hand, combining the main features of GMRES-DR and GCRO, in 2006, an algorithm called GCRO-DR was derived in~\cite{Parks2006} for sequences of linear systems.
When solving a single linear system, GCRO-DR is mathematically equivalent to GMRES-DR.

It was Gutknecht~\cite{Gutknecht2012} who attempted to present deflated and augmented Krylov subspace methods in a common framework.
Also, Gutknecht investigated the breakdown conditions of deflated GMRES; we also refer the interested reader to the companion paper~\cite{Gaul2013}.
These two papers might provide good starting point for further investigations of deflation and augmentation.

\subsection{Weighted inner products}

In 1998, it was shown by Essai~\cite{Essai1998} that a suitably chosen inner product could improve the performance of restarted MGS-GMRES.
Let~$D=\diag(\delta_1,\dots,\delta_N)$ with $\delta_i>0$.
The $D$-inner product is defined as $(u,v)_D=v^\intercal Du$ for all $u,v\in\mathbb{R}^N$, and then the associated $D$-norm is $\norm{u}_D=\sqrt{u^\intercal Du}$.
From an algorithmic viewpoint, for example in Algorithm~\ref{alg:arnoldi:mgs}, one only needs to replace lines 5 and 7 with the new inner product and norm.
The resulting algorithm shall be denoted by WGMRES.
In~\cite{Essai1998}, the diagonal matrix~$D$ was chosen as~$\delta_i=\sqrt{N}\abs{r_{0,i}}/\norm{r_0}$ and updated at each restart.

Preconditioned variants were studied in~\cite{Cao2004}.
Later, a modified version motivated by the idea of~\cite{Ayachour2003} was proposed in~\cite{SNajafi2008}, in which $\delta_i$ are chosen randomly and there is no need to use Givens rotations.
In~\cite{Niu2010}, WGMRES was combined with the augmentation technique, and then the new variant was compared with LGMRES as illustrated in Algorithm~\ref{alg:lgmres}.
On the other hand, Embree et al.~\cite{Embree2017} in 2017 applied weighting to GMRES-DR and observed that the mixed scheme could be better than using deflation alone or weighting alone.
In 2014, G{\"u}ttel and Pestana~\cite{Guttel2014} proposed a new variant that applied the orthogonalization process to the transformed matrix~$D^{\frac{1}{2}}AD^{-\frac{1}{2}}$ and starting vector~$D^{\frac{1}{2}}r_0$ with respect to the Euclidean inner product, which is mathematically equivalent to WGMRES.
For more insight in nonstandard inner products, see~\cite{Pestana2013,Guttel2014,Embree2017} and references therein.

\subsection{Inexact matrix-vector products}

An inexact process might be beneficial if the computation of the matrix-vector product~$Av_j$ is time-consuming.
From a mathematical point of view, an inexact matrix-vector product can be modeled by~$(A+E_j)v_j$ where $E_j$ is a perturbation matrix at~$j$th iteration.
From~\eqref{eq:arnoldi}, one obtains that
\[
AV_n + [E_1v_1,\dots,E_nv_n] = V_{n+1}\overbar{H}_n.
\]
The inexact GMRES algorithm was first proposed in a technical report by Bouras and Frayss\'e~\cite{Bouras2000} in~2000 (later published in~\cite{Bouras2005} in~2005), and then analyzed by Simoncini and Szyld~\cite{Simoncini2003,Simoncini2005}, van den Eshof and Sleijpen~\cite{vandenEshof2004}, Giraud et al.~\cite{Giraud2007}, and Sifuentes et al.~\cite{Sifuentes2013}.
A relaxation strategy on the inner accuracy of GMRES can be found in~\cite{Bouras2005}.
It was shown in~\cite{Simoncini2003} that the early matrix-vector products must be computed with high accuracy, but can be relaxed as the iteration progresses, which explains the term ``relaxation''.
An interesting review of inexact GMRES was given in~\cite{Simoncini2007}; see~\cite{Sidje2011} for performance evaluation.
In 2013, Dolgov~\cite{Dolgov2013} employed relaxation techniques to improve the performance of tensor GMRES; see~\cite{Ballani2013,Dolgov2013} for more details on tensor formats.

\subsection{Mixed precision}

On modern computer architectures, single precision arithmetic (32 bits) is cheaper than double precision arithmetic (64 bits) in the sense that communication and computation costs grow with the size of the floating point format~\cite{Abdelfattah2021}.
In the 2008 revision of the IEEE standard, half precision (16 bits) was defined as a storage format.
Mixed-precision numerical algorithms can date back to more than five decades ago~\cite{Moler1967}.
Recent advances in hardware-level support (e.g.,~NVIDIA V100 GPU) for low-precision arithmetic and the increasing impact of communication costs have promoted the reuse of mixed-precision techniques; see~\cite{Abdelfattah2021} for a recent survey.

It is very natural to extend Algorithm~\ref{alg:gmres:mgs} to the case of mixed precision.
Here, we only consider the two-precision algorithm, which exploits single precision everywhere, except that
\begin{enumerate}
\item The original system $Ax=b$ must be stored in double precision;
\item Although $V_ny$ could be computed in single precision, the update of~$x_n$ should be done in double precision;
\item For the restarted version, the computation of residual vector should also be done in double precision.
\end{enumerate}
A recent paper by Gratton et al.~\cite{GrattonX1} showed that using low precision in MGS-GMRES could achieve the same convergence rate and final accuracy as the full precision variant.
They provided a theoretical analysis for the non-restarted algorithm by means of inexact inner products.
Then, Lindquist et al.~\cite{Lindquist2020} focused on the experimental aspects of restarted mixed-precision GMRES algorithms, including MGS-GMRES and a variant based on CGS with one reorthogonalization (CGS2); see~\cite{Giraud2005,Giraud2005b} for error analysis of CGS2.
\begin{algorithm}
\caption{CGS2-Arnoldi\label{alg:arnoldi:cgs2}}
\begin{algorithmic}[1]
  \State $v_1=r_0/\norm{r_0}$
  \For {$j=1,2,\dots,n$}
    \State $w_j=Av_j$
    \State $\hat{h}=V_j^\mathsf{T}w_j,\quad w_j=w_j-V_j\hat{h}$
    \State $\tilde{h}=V_j^\mathsf{T}w_j,\quad w_j=w_j-V_j\tilde{h}$
    \State $[h_{1,j},\dots,h_{j,j}]^\mathsf{T}=\hat{h}+\tilde{h}$
    \State $h_{j+1,j}=\norm{w_j}$,\quad $v_{j+1}=w_j/h_{j+1,j}$
  \EndFor
\end{algorithmic}
\end{algorithm}
In Algorithm~\ref{alg:arnoldi:cgs2}, we illustrate the CGS2-Arnoldi procedure.
Motivated by the idea in~\cite{GrattonX1}, a theoretical result for CGS2-GMRES was provided in~\cite{Lindquist2021}, along with more experimental results.

Another interesting class of mixed-precision algorithms involves iterative refinement where low-precision GMRES is employed as inner solver.
In~1992, Turner and Walker~\cite{Turner1992} showed that GMRES($m$) can be regarded as an iterative refinement process.
In~2009, Arioli and Duff~\cite{Arioli2009} proved that the FGMRES algorithm preconditioned by low-precision factorization is backward stable.
\begin{algorithm}
\caption{GMRES-IR\label{alg:gmres-ir}}
\begin{algorithmic}[1]
  \State Compute the LU factorization $A=LU$ in single precision
  \State Solve $Mx_0=b$ in single precision with $M=LU$
  \State $r_0=b-Ax_0$ in double precision
  \For {$k=0,1,\dots$}
    \State Solve $M^{-1}Ad_k=M^{-1}r_k$ by a GMRES algorithm in single precision with $M=LU$
    \State $x_{k+1}=x_k+d_k$ in double precision
    \State $r_{k+1}=b-Ax_{k+1}$ in double precision
    \State Check convergence
  \EndFor
\end{algorithmic}
\end{algorithm}
We sketch GMRES-based iterative refinement (GMRES-IR) in Algorithm~\ref{alg:gmres-ir}; see,~e.g.,~\cite{Carson2018c} for more details.
Note that this is just a simple example illustrated using two precisions, for which three or more precisions could also be exploited by relaxing some of the operations.
Carson and Higham~\cite{Carson2017} in 2017 gave a new forward analysis of iterative refinement and suggested the use of GMRES-IR.
Then, in~\cite{Carson2018c}, they introduced and analyzed an iterative refinement scheme that can possibly exploit three different precisions.
It was shown that GMRES-IR could provide the required relative accuracy, even though the LU factorization in single precision might be of low quality.
The error analysis in~\cite{Carson2017,Carson2018c} was recently generalized in~\cite{AmestoyX1} to the case where the GMRES solver and matrix-vector multiplications in line~5 may proceed in independent precisions, yielding a possibility of using five precisions.
On the other hand, the overflow and underflow problems associated with mixed-precision algorithms were investigated by Higham et al.~\cite{Higham2019}.

The survey article~\cite{Abdelfattah2021} provides additional details and references on mixed-precision techniques.
Moreover, although not directly connected to GMRES, the studies in~\cite{Yamazaki2015} and~\cite{Yang2021} can deepen our insight into mixed-precision orthogonalization algorithms.

\section{Parallel algorithms}\label{sec:pa}

On modern computer architectures, communication is expensive relative to computation.
Sparse matrix-vector products (SpMVs) and dot products, especially the latter, often limit the possible speedups for parallel algorithms.
Therefore, the high cost of global reduction operations in the orthogonalization process makes GMRES hard to parallelize.
To alleviate performance bottleneck, much work has focused on ways to reduce or hide communication costs.
Early work concerning the parallelization of GMRES can be found in,~e.g.,~\cite{Walker1988,diBrozolo1989,Shadid1992,daCunha1994}; see also~\cite{Demmel1993} for a general discussion on parallel algorithms.
Here we focused on the $s$-step GMRES algorithm, and then present the idea of pipelining.
These approaches have gained popularity in the last decade.
Finally, some other potential strategies are included at the end of this section.

\subsection{Matrix powers and polynomial basis}

In~\cite{Demmel2007} and some subsequent publications~\cite{Demmel2008,Mohiyuddin2009,Hoemmen2010,Ballard2014}, a Krylov subspace algorithm is thought to contain several ``kernels'', each kernel being an important and time-consuming building block of that algorithm.
In GMRES, matrix powers computation is thus considered as a kernel, which computes $s$ vectors~$Av,\dots,A^sv$.
The most trivial implementation requires $s$ messages between each processor and its neighbors, whereas a well-designed matrix powers kernel (MPK) may require only $1+o(1)$ message at the cost of redundant computation and storage consumption; see~\cite{Demmel2007}, in which data dependencies and performance models are discussed in terms of both distributed and shared memory systems.
Nonetheless, this technique has several disadvantages.
First, MPK introduces additional overheads.
Second, it's still difficult to apply effective preconditioners with MPK~\cite{Hoemmen2010,Yamazaki2017}.
Third, basis vectors will become more and more linearly dependent~\cite{deSturler1991}.

The last drawback can be partially corrected by using polynomial basis, that is,
\[
W_s = [w_1,\dots,w_s] = [\phi_0(A)w,\phi_1(A)w,\dots,\phi_{s-1}(A)w],
\]
where $\phi_j(z)$ is a polynomial of degree~$j$.
Choosing
\begin{equation}\label{eq:basis:m}
\phi_j(z)=z^j
\end{equation}
yields the monomial basis, which converges to the eigenvector corresponding to the dominant eigenvalue of~$A$, and thus, as mentioned above, will rapidly become numerically rank deficient.
Two alternatives have been suggested.
One is to use Newton polynomials
\begin{equation}\label{eq:basis:n}
\phi_0(z)=1,\quad \phi_j(z)=\varrho_j(z-\theta_j)\phi_{j-1}(z),
\end{equation}
as suggested by Bai et al.~\cite{Bai1994} in 1994, in which an $s$-step GMRES algorithm was proposed.
Here, $\varrho_j$ are scaling factors, which are generally taken as~$1$, and $\theta_j$ are eigenvalue estimates computed by a few iterations of the Arnoldi process.
The other is Chebyshev polynomials, which can be defined by
\begin{equation}\label{eq:basis:c}
\phi_0(z)=1,\quad \phi_1(z)=\frac{1}{2\gamma}(z-\zeta)\phi_0(z),\quad \phi_j(z)=\frac{1}{\gamma}((z-\zeta)\phi_{j-1}(z)-\frac{\tau^2}{4\gamma}\phi_{j-2}(z)),
\end{equation}
where the spectrum is assumed to be circumscribed by the rectangle~$\{z=z_1+\iota z_2:\abs{z_1-\zeta}\le\xi_1,\abs{z_2}\le\xi_2\}$ with $\iota$ the imaginary unit.
Then, we choose $\gamma=\max\{\xi_1,\xi_2\}$ and $\tau=\sqrt{\xi_1^2-\xi_2^2}$, so that $\zeta\pm\tau$ are the foci of the ellipse for the spectrum.
Chebyshev polynomials for GMRES were discussed by Joubert and Carey~\cite{Joubert1992,Joubert1992b} in 1992; more details of Krylov subspace bases can be found in~\cite{Hoemmen2010,Philippe2012,Carson2015b}.

\subsection{$s$-step GMRES}

The idea of performing $s$ iterations at once can date back to the early 1950s, as quoted in~\cite{Forsythe1968}.
In 1990, Chronopoulos and Kim~\cite{Chronopoulos1990} developed an $s$-step GMRES algorithm using monomial basis.
A Newton basis $s$-step GMRES algorithm was considered by Bai et al.~\cite{Bai1994}, in which modified Leja ordering and basis scaling were employed to enhance robustness.
They also suggested the use of parallel QR factorization, and this was implemented by Erhel~\cite{Erhel1995} in 1995 by means of the so-called RODDEC algorithm.
On the other hand, Joubert and Carey~\cite{Joubert1992,Joubert1992b} in 1992 considered an $s$-step GMRES algorithm using Chebyshev basis along with an MPK for 2D regular meshes; see also~\cite{deSturler1995}.
In 2009, Mohiyuddin et al.~\cite{Mohiyuddin2009} introduced an improved variant, for which a new kernel called tall skinny QR was derived to achieve the optimal reduction of communication while maintaining the accuracy; see~\cite{Demmel2012}.
Then, the thesis by Hoemmen~\cite{Hoemmen2010} gave a complete treatise of $s$-step GMRES.

Let $\overbar{B}_s$ be an~$(s+1)\times s$ matrix such that
\begin{equation}\label{eq:bc}
AW_s=W_{s+1}\overbar{B}_s.
\end{equation}
It is called a basis conversion matrix, by which the three bases with respect to~\eqref{eq:basis:m}--\eqref{eq:basis:c} can be exploited in an elegant manner.
It was shown in~\cite{Hoemmen2010} that the QR factorization of~$W_s$ could produce the same orthogonal matrix as that in~\eqref{eq:arnoldi}.
Then, it follows that
\begin{equation}\label{eq:wsqr}
W_s = V_sT_s,
\end{equation}
where $T_s\in\mathbb{R}^{s\times s}$ is an upper triangular matrix.
Along with~\eqref{eq:arnoldi:1} and~\eqref{eq:bc}, this implies
\[
\overbar{H}_s = T_{s+1}\overbar{B}_sT_s^{-1}.
\]
The above formulation can potentially bring two benefits.
First, one could compute $W_{s+1}$ with less communication than in the classical case.
Second, the QR factorization $W_{s+1}=V_{s+1}T_{s+1}$ could be performed in parallel.

The tall skinny QR (TSQR) algorithm studied in~\cite{Hoemmen2010,Demmel2012} is a good choice for reducing communication when there are many more rows than columns.
Let $W_s$ be vertically partitioned into~$\nu_0$ submatrices such that the number of rows of each submatrix is larger than~$s$.
Then, we perform the QR factorization of each block entry, that is,
\[
W_s = \left(\begin{array}{c} W_{s,1} \\ W_{s,2} \\ \vdots \\ W_{s,\nu_0} \end{array}\right)
= \left(\begin{array}{c} Q_1^{(0)}R_1^{(0)} \\ Q_2^{(0)}R_2^{(0)} \\ \vdots \\ Q_{\nu_0}^{(0)}R_{\nu_0}^{(0)} \end{array}\right)
= \left(\begin{array}{cccc} Q_1^{(0)} & & & \\ & Q_2^{(0)} & & \\ & & \ddots & \\ & & & Q_{\nu_0}^{(0)} \end{array}\right)
\left(\begin{array}{c} R_1^{(0)} \\ R_2^{(0)} \\ \vdots \\ R_{\nu_0}^{(0)} \end{array}\right).
\]
Once such initialization is complete, in each cycle we perform the QR factorization of vertically stacked upper triangular submatrices
\[
\left(\begin{array}{c} R_1^{(0)} \\ R_2^{(0)} \end{array}\right) = Q_1^{(1)}R_1^{(1)},\quad \left(\begin{array}{c} R_3^{(0)} \\ R_4^{(0)} \end{array}\right) = Q_2^{(1)}R_2^{(1)},\quad \left(\begin{array}{c} R_1^{(1)} \\ R_2^{(1)} \end{array}\right) = Q_1^{(2)}R_1^{(2)},\quad\dots,
\]
with $R_i^{(j)}\in\mathbb{R}^{s\times s}$, and define
\[
  Q^{(j)} = \left(\begin{array}{cccc} Q_1^{(j)} & & & \\ & Q_2^{(j)} & & \\ & & \ddots & \\ & & & Q_{\nu_j}^{(j)} \end{array}\right),\quad
  R^{(j)} = \left(\begin{array}{c} R_1^{(j)} \\ R_2^{(j)} \\ \vdots \\ R_{\nu_j}^{(j)} \end{array}\right).
\]
Here, $\nu_j=\lceil\nu_{j-1}/2\rceil$ for~$j=1,2,\dots,\lceil\log_2{\nu_0}\rceil$.
If~$\nu_{j-1} \bmod 2\ne 0$, then~$Q_{\nu_j}^{(j)}=I$ and~$R_{\nu_j}^{(j)}=R_{\nu_{j-1}}^{(j-1)}$.
Algorithm~\ref{alg:tsqr} shows this procedure.
\begin{algorithm}
\caption{TSQR\label{alg:tsqr}}
\begin{algorithmic}[1]
  \For {$i=1,2,\dots,\nu_0$}
    \State Compute QR factorization~$W_{s,i}=Q_i^{(0)}R_i^{(0)}$
  \EndFor
  \For {$j=1,2,\dots,\lceil\log_2{\nu_0}\rceil$}
    \State $\hat{\nu}_j=\lfloor\nu_{j-1}/2\rfloor,\quad \nu_j=\lceil\nu_{j-1}/2\rceil$
    \For {$i=1,2,\dots,\hat{\nu}_j$}
      \State Compute QR factorization~$\left(\begin{array}{c}R_{2i-1}^{(j-1)} \\ R_{2i}^{(j-1)}\end{array}\right)=Q_i^{(j)}R_i^{(j)}$
    \EndFor
    \If {$\nu_j\ne\hat{\nu}_j$}
      \State $Q_{\nu_j}^{(j)}=I,\quad R_{\nu_j}^{(j)}=R_{\nu_{j-1}}^{(j-1)}$
    \EndIf
  \EndFor
\end{algorithmic}
\end{algorithm}
Notice that~$Q^{(\lceil\log_2{\nu_0}\rceil)}=Q_1^{(\lceil\log_2{\nu_0}\rceil)}$ and~$R^{(\lceil\log_2{\nu_0}\rceil)}=R_1^{(\lceil\log_2{\nu_0}\rceil)}$.
Then, along with~\eqref{eq:wsqr}, it follows that
\[
V_s = Q^{(0)}Q^{(1)}\dots Q^{(\lceil\log_2{\nu_0}\rceil)},\quad T_s = R^{(\lceil\log_2{\nu_0}\rceil)},
\]
yielding the desired factorization of~$W_s$.

A large $s$ may lead to stability problems even when a polynomial basis is used.
In~\cite{Mohiyuddin2009,Hoemmen2010} it was thus suggested to combine TSQR with the block Gram-Schmidt (BGS) algorithm, such that $s$ can be chosen independently of the restart length~$m$.
As a result, one can pursue both numerical stability and convergence rate.
Now, assume that $m=s\cdot t$.
The BGS process can be written as
\[
\mathfrak{R}_s=\mathfrak{V}_{s(j-1)}^\mathsf{T} W_s,\quad W_s=W_s-\mathfrak{V}_{s(j-1)}\mathfrak{R}_s,
\]
where $\mathfrak{V}_{s(j-1)}$ denotes the orthogonal matrix that contains all previous orthogonal basis vectors established in this cycle.
After performing the TSQR procedure, we can get~$\mathfrak{V}_{sj}=[\mathfrak{V}_{s(j-1)}, V_s]$.
In addition, let $\mathfrak{H}_{sj}$ be the upper Hessenberg matrix corresponding to~$\mathfrak{V}_{sj}$.
Then the BGS-based $s$-step GMRES satisfies the relation
\[
A\mathfrak{V}_{sj} = \mathfrak{V}_{sj+1}\overbar{\mathfrak{H}}_{sj},
\]
which is equivalent to~\eqref{eq:arnoldi:1} if choosing~$n=sj$.
Here, BGS amounts to a block version of CGS.
For more details on BGS, we refer the reader to~\cite{Carson2021,Carson2022} and references therein.
BGS and TSQR are both regarded as kernels (see~\cite{Mohiyuddin2009,Hoemmen2010}).

Notice that the last basis vector of~$\mathfrak{V}_{sj+1}$ can be used as the starting vector for the next~$W_s$.
Hence, we use the acute accent to indicate a shift of indices such as $\acute{W}_s=[w_2,w_3,\dots,w_{s+1}]$, on which BGS and TSQR would be performed.
\begin{algorithm}
\caption{$s$-step GMRES\label{alg:s-gmres}}
\begin{algorithmic}[1]
  \State $\beta=\norm{r_0},\quad w_1=r_0/\beta$
  \For {$j=1,2,\dots,t$}
    \State Select~$\overbar{B}_s$ and compute $\acute{W}_s=[w_2,w_3,\dots,w_{s+1}]$
    \If {$j=1$}
      \State Compute QR factorization~$W_{s+1}=V_{s+1}T_{s+1}$
      \State $\overbar{H}_s=T_{s+1}\overbar{B}_sT_s^{-1}$
      \State $\mathfrak{V}_{s+1}=V_{s+1},\quad \overbar{\mathfrak{H}}_s=\overbar{H}_s,\quad \overbar{\mathfrak{B}}_s=\overbar{B}_s$
    \Else
      \State $\acute{\mathfrak{R}}_s=\mathfrak{V}_{s(j-1)+1}^\mathsf{T} \acute{W}_s,\quad \acute{W}_s=\acute{W}_s-\mathfrak{V}_{s(j-1)+1}\acute{\mathfrak{R}}_s$
      \State Compute QR factorization~$\acute{W}_s=\acute{V}_s\acute{T}_s$
      \State $\overbar{\mathfrak{B}}_{sj}=\left(\begin{array}{cc}\mathfrak{H}_{s(j-1)} & 0 \\ \eta_{j-1}e_1e_{s(j-1)}^\mathsf{T} & \overbar{B}_s\end{array}\right),\quad \mathfrak{T}_{sj+1}=\left(\begin{array}{cc}I_{s(j-1)+1} & \acute{\mathfrak{R}}_s \\ 0 & \acute{T}_s\end{array}\right)$
      \State $\overbar{\mathfrak{H}}_{sj}=\mathfrak{T}_{sj+1}\overbar{\mathfrak{B}}_{sj}\mathfrak{T}_{sj}^{-1},\quad \mathfrak{V}_{sj+1}=[\mathfrak{V}_{s(j-1)+1}, \acute{V}_s]$
    \EndIf
  \EndFor
  \State Compute~$y$ such that $\norm{\beta e_1-\overbar{\mathfrak{H}}_{st} y}$ is minimized
  \State $x_{st}=x_0+\mathfrak{V}_{st} y$
\end{algorithmic}
\end{algorithm}
Algorithm~\ref{alg:s-gmres} illustrates the $s$-step GMRES algorithm in a restart cycle, which is mathematically equivalent to Algorithm~\ref{alg:gmres:mgs}.
Indeed, the lack of consensus on the meaning of ``$s$-step GMRES'' often leads to confusion.
In~\cite{Mohiyuddin2009,Hoemmen2010} the new formulation with MPK, BGS and TSQR was called communication-avoiding GMRES (CA-GMRES).
In our context, Algorithm~\ref{alg:s-gmres} depicts the $s$-step GMRES algorithm, whereas CA-GMRES refers to the MPK-based $s$-step GMRES algorithm, namely, a matrix powers kernel is employed in line~3, where $s$ basis vectors can be computed for the same communication cost as computing 1 basis vector.
Here $\mathfrak{T}_{sj}$ is the $sj\times sj$ principal submatrix of~$\mathfrak{T}_{sj+1}$.
In line~11, $\eta_{j-1}$ denotes the bottommost value of the $s(j-1)$th column of~$\overbar{\mathfrak{H}}_{s(j-1)}$.
For the sake of clarity, we denote by~$I_s$ the identity matrix of size~$s$.
Further, the upper Hessenberg matrix~$\overbar{\mathfrak{H}}_{sj}$ in line~12 can be updated in a more efficient way; see~\cite{Hoemmen2010} for more details.

A full account of TSQR can be found in~\cite{Demmel2012}, while Hoemmen's thesis~\cite{Hoemmen2010} contains further discussion of QR algorithms.
Parallel version of deflated GMRES was studied by Wakam and Erhel~\cite{Wakam2013} in 2013, who combined Newton basis $s$-step GMRES with adaptive augmentation.
Similar topics were discussed independently in~\cite{Yamazaki2014c}; we refer the reader to~\cite{Wakam2013} and~\cite{Yamazaki2014c} for algorithms and details.
In 2017, Imberti and Erhel~\cite{Imberti2017} proposed an $s$-step GMRES algorithm where the block size~$s$ is variable.
They suggested to use a predetermined sequence based on the Fibonacci numbers.
Experiments on multiple GPUs can be found in~\cite{Yamazaki2014b}, with a particular focus on orthogonalization strategies; see also the experiments in~\cite{Yamazaki2017,Yamazaki2020}.
Moreover, polynomial preconditioning could be used with $s$-step GMRES; see~\cite{Loe2020} for more discussion.

\subsection{Pipelined GMRES}

In 2013, Ghysels et al.~\cite{Ghysels2013} proposed a pipelined variant of GMRES which only requires a single nonblocking reduction per iteration and can interleave dot products with SpMVs to hide communication latency.
An important observation is that \eqref{eq:cgs3} can be rewritten via the Pythagorean theorem as
\begin{equation}\label{eq:cgs-p}
h_{j+1,j} = \sqrt{\norm{Av_j}^2 - \sum_{i=1}^j h_{i,j}^2}.
\end{equation}
The CGS process involving~\eqref{eq:cgs-p} could be denoted by CGS-P (see,~e.g.,~\cite{Smoktunowicz2006,Carson2021}), in which one finds that only one reduction is required, whereas two are required in the traditional CGS-Arnoldi algorithm.
This results in the single-reduction GMRES algorithm.
At the $j$th iteration, the orthogonalization process can be summarized as follows:
\begin{enumerate}
\item $w_j=Av_j$;
\item Compute $h_{i,j}=(w_j,v_i)$ for $i=1,2,\dots,j$ and $(w_j,w_j)$ by a global reduction;
\item $h_{j+1,j}=\sqrt{(w_j,w_j)-\sum_{i=1}^jh_{i,j}^2}$;
\item $v_{j+1}=(w_j-\sum_{i=1}^jh_{i,j}v_i)/h_{j+1,j}$.
\end{enumerate}
In practice for stability reasons, a polynomial basis should be used, and a reorthogonalization process could be helpful when a breakdown occurs in step~3 (see~\cite{Ghysels2013}).

It was observed that $w_{j+1}$ can be obtained by means of an updating formula.
Along with a subtle rearrangement of operations, we obtain a communication-hiding variant shown in Algorithm~\ref{alg:p-gmres}.
\begin{algorithm}
\caption{Pipelined GMRES\label{alg:p-gmres}}
\begin{algorithmic}[1]
  \State $\beta=\norm{r_0},\quad v_1=r_0/\beta,\quad w_1=(A-\theta I)v_1$
  \For {$j=1,2,\dots,n$}
    \State Compute $h_{i,j}=(w_j,v_i)$ for $i=1,2,\dots,j$ and $\varsigma=(w_j,w_j)$ using nonblocking reductions
    \State $u=Aw_j$
    \State $h_{j+1,j}=\sqrt{\varsigma-\sum_{i=1}^jh_{i,j}^2}$
    \State $v_{j+1}=(w_j-\sum_{i=1}^jh_{i,j}v_i)/h_{j+1,j}$
    \State $h_{i,j}=h_{i,j}+\theta$
    \State $w_{j+1}=(u-\sum_{i=1}^jh_{i,j}w_i)/h_{j+1,j}$
  \EndFor
  \State Compute~$y$ such that $\norm{\beta e_1-\overbar{H}_n y}$ is minimized
  \State $x_n=x_0+V_ny$
\end{algorithmic}
\end{algorithm}
The single-stage pipelined GMRES algorithm, here loosely called pipelined GMRES, can achieve a speedup factor of at most~2.
We can see in the first line that the Newton basis is used to improve numerical stability.
In the for-loop, Newton polynomials are implicitly involved.
If a dot product entails higher latency than an SpMV, then a greater speedup might be possible by overlapping a nonblocking reduction with at most $l$ SpMVs, in other words, using deeper pipelining to hide a dot-product latency of at most $l$ iterations; see~\cite{Ghysels2013} for the algorithm and a detailed derivation.
The $l$-stage variant, however, does not decrease the number of global reductions and can introduce more source of instability.

Note that pipelining allows us to easily use preconditioners, whereas $s$-step algorithms would encounter some difficulties (see,~e.g,~\cite{Ghysels2013,Yamazaki2017}).
In 2016, Sanan et al.~\cite{Sanan2016} presented a pipelined variant of FGMRES, which was also mentioned in~\cite{Ghysels2013}, and expressed doubts about the combination of deeper pipelining with nonlinear preconditioners.
The idea of combining the strengths of the pipelined and $s$-step algorithms was considered by Yamazaki et al.~\cite{Yamazaki2017}.
It was argued in~\cite{Loe2020} that polynomial preconditioners would be well suited to pipelined GMRES.
Besides, pipelining was used in~\cite{Morgan2021} to study the performance variability on high-performance computing systems.
We also refer the reader to~\cite{Eller2019} for a detailed account of nonblocking Krylov solvers.

\subsection{Low-sync GMRES}

In 2007, Hern{\'a}ndez et al.~\cite{Hernandez2007} developed three variants of the Arnoldi process for distributed-memory platforms.
In particular, the so-called Arnoldi with delayed reorthogonalization (ADR) algorithm could reduce synchronization points by rearranging and deferring operations.
Numerical experiments reported in~\cite{Hernandez2007} show that ADR provides a large performance improvement when the number of processors was sufficiently increased, but suffers from numerical instability; see~\cite{Hoemmen2010} for a comparison of ADR and $s$-step strategies.

Recently, Swirydowicz et al.~\cite{Swirydowicz2021} suggested a new variant of MGS-GMRES that requires only one reduction per iteration.
The main idea proposed is to apply the delayed normalization technique of ADR to the inverse compact WY representation of MGS.
Recall that MGS can be viewed as a composition of rank-one projectors.
Swirydowicz et al.~\cite{Swirydowicz2021} noted the inverse compact WY form of MGS (ICWY-MGS), which was partly anticipated by Bj{\"o}rck's~\cite{Bjorck1967} more than five decades ago, but laid dormant for a long time.
Defining
\[
L_n = \left(\begin{array}{cc}
L_{n-1} & 0 \\
(V_{n-1}^\mathsf{T}v_n)^\mathsf{T} & 0
\end{array}\right),\quad L_1=0,
\]
it follows that
\[
\hat{P}_n = (I-v_nv_n^\intercal)\dots(I-v_1v_1^\intercal)=I-V_n\hat{T}_nV_n^\mathsf{T},\quad \hat{T}_n=(I+L_n)^{-1},
\]
where $L_n$~is a strictly lower triangular matrix in finite precision arithmetic; see~\cite{Swirydowicz2021}.
Only one reduction is needed for ICWY-MGS to perform the projection.
In general, line~7 in Algorithm~\ref{alg:arnoldi:mgs} requires an extra global reduction for the normalization.
Now, let~$l_{i,j}$ and~$h_{i,j}$ denote respectively the $(i,j)$ entries of~$L_n$ and~$\overbar{H}_n$.
The delayed normalization technique described in~\cite{Hernandez2007} was employed in~\cite{Swirydowicz2021}, resulting in a low-sync GMRES algorithm as shown in Algorithm~\ref{alg:gmres:icwy}.
\begin{algorithm}
\caption{Low-sync GMRES\label{alg:gmres:icwy}}
\begin{algorithmic}[1]
  \State $\beta=\norm{r_0},\quad v_1=r_0/\beta,\quad h_{1,1}=(w_1,v_1),\quad w_1=w_1-h_{1,1}v_1$
  \For {$j=2,3,\dots,n+1$}
    \State $w_j=Aw_{j-1}$
    \State $l_{j,1:j-1}=(V_{j-1}^\mathsf{T}w_{j-1})^\mathsf{T},\quad u=[V_{j-1},w_{j-1}]^\mathsf{T}w_j,\quad h_{j,j-1}=\norm{w_{j-1}}$
    \State $v_j=w_{j-1}/h_{j,j-1}$
    \State $l_{j,1:j-1}=l_{j,1:j-1}/h_{j,j-1},\quad u=u/h_{j,j-1},\quad u_j=u_j/h_{j,j-1},\quad w_j=w_j/h_{j,j-1}$
    \State $h_{1:j,j}=(I+L_j)^{-1}u$
    \State $w_j=w_j-V_jh_{1:j,j}$
  \EndFor
  \State Compute~$y$ such that $\norm{\beta e_1-\overbar{H}_n y}$ is minimized
  \State $x_n=x_0+V_ny$
\end{algorithmic}
\end{algorithm}
We can see that delaying the normalization of~$w_{j-1}$ brings additional scaling (line~6), but can merge the normalization with dot products (line~4).
As a result, only one global reduction per iteration is required.

Several formulations of Gram-Schmidt were studied in~\cite{Swirydowicz2021} and~\cite{BielichX1} as they can minimize the number of global reductions.
For example, a single-reduction CGS2 variant can also be derived by means of delayed normalization.
Numerical results reported promising performance of the low-sync Gram-Schmidt and GMRES algorithms.
A good survey of Gram-Schmidt up to 2010 was given by Leon et al.~\cite{Leon2013} in 2013; we also refer the interested reader to~\cite{Barlow2013,Barlow2019,Carson2021,Carson2022} and their references for the discussion of block Gram-Schmidt algorithms.
On the other hand, a low-sync GMRES variant based on Householder QR was derived by Walker~\cite{Walker1988} in 1988.
The single-reduction orthogonalization procedures were applied to $s$-step and pipelined algorithms by Yamazaki et al.~\cite{Yamazaki2020}, together with numerical experiments on their performance.
The recent survey article~\cite{Abdelfattah2021} that were intended to provide a road map for mixed-precision algorithms contains also a perspective on mixed-precision ICWY-MGS.

\section{Extensions to other problems}\label{sec:ext}

This section provides a concise literature review on GMRES solvers for other related problems.
First, a more complicated but interesting problem is to solve systems with multiple right-hand sides
\begin{equation}\label{eq:mrhs}
AX=B,
\end{equation}
where the $N\times S$ matrix~$B$ is a collection of right-hand sides with~$S\ll N$.
It is known that block Krylov algorithms can handle all columns at once, and thus has a few practical advantages in terms of data locality and commnication efficiency.
Block GMRES was described by Vital~\cite{Vital1990} in 1990 and then studied by Simoncini and Gallopoulos~\cite{Simoncini1995,Simoncini1996}; see also~\cite{Gu2001,Saad2003,Morgan2005,Robbe2006,Gutknecht2009,Calandra2012,Calandra2013,Agullo2014,ADaas2019,Frommer2020,Kubinova2020,Tajaddini2021}.
Block GMRES has also been extended to handle the single right-hand side system~\eqref{eq:ls}; see,~e.g.,~\cite{Chapman1997,Baker2006,Greif2017}.
As an alternative, Jbilou et al.~\cite{Jbilou1999} in 1999 proposed a global GMRES algorithm for solving~\eqref{eq:mrhs}; see also~\cite{Bouyouli2006,Elbouyahyaoui2009}.

In 2006, Parks et al.~\cite{Parks2006} proposed a recycled GMRES algorithm, called GCRO-DR, for solving a sequence of linear systems
\[
A^{(i)}x^{(i)} = b^{(i)},
\]
where the matrix changes slowly from one system to the next, but the cumulative change over many iterations can be large.
On the other hand, Frommer and Gl{\"a}ssner~\cite{Frommer1998b} in 1998 considered the shifted coefficient matrix~$A+\sigma I$ and proposed a shifted GMRES algorithm for this case.
The solution of sequences of shifted linear systems
\begin{equation}\label{eq:sls}
(A+\sigma_i I)x^{(i)} = b
\end{equation}
was studied in~\cite{Soodhalter2014,Soodhalter2016,Soodhalter2016b} by combining shifted GMRES with subspace recycling.
We refer to the survey by Soodhalter et al.~\cite{Soodhalter2020} for more details on recycling strategies.
A multi-preconditioned GMRES algorithm was proposed in~\cite{Greif2017}, which when combined with the flexible approach developed in~\cite{Saibaba2013}, could be applied to~\eqref{eq:sls}; see~\cite{Bakhos2017}.
GMRES could have broader applications for more general shifted systems; see~\cite{Darnell2008,Soodhalter2014,Soodhalter2016,Sun2018,Frommer2020} for more details.
Now consider the Sylvester equation
\begin{equation}\label{eq:syl}
AX+XD=B,
\end{equation}
where $A\in\mathbb{R}^{N\times N}$, $D\in\mathbb{R}^{S\times S}$, and $B\in\mathbb{R}^{N\times S}$.
A block GMRES algorithm for solving~\eqref{eq:syl} was described by Robb{\'e} and Sadkane~\cite{Robbe2002} in 2002; see also~\cite{Simoncini2016,Zadeh2019}.
It is known that a family of shifted systems can be interpreted as a Sylvester equation of the form~\eqref{eq:syl} where $D$ is a diagonal matrix; see~\cite{Soodhalter2016,Elbouyahyaoui2021}.
GMRES has also been extended to solve least square problems~\cite{Hayami2010}, tensor equations~\cite{Chen2012}, and quaternion linear systems~\cite{Jia2021}.

Going back to the original system~\eqref{eq:ls}, we point the reader to~\cite{Fischer1998,Benzi2005,Frommer2020b,Southworth2020}, where the coefficient matrix~$A$ has a $2\times2$ block structure.
There are practical cases where the matrix~$A$ is singular or nearly singular.
The article~\cite{Brown1997} by Brown and Walker in 1997 seems to be the first to study the behavior of GMRES for solving singular or ill-conditioned systems.
Later, the Drazin-inverse solution of singular linear systems was discussed by Ipsen and Meyer~\cite{Ipsen1998} and Sidi~\cite{Sidi2001}; see also~\cite{Calvetti2002,Reichel2005,Elden2012,Greenbaum2018,Morikuni2018,Gazzola2019}.
In parallel to these continuing efforts, a range restricted GMRES algorithm geared again toward singular systems was developed by Calvetti et al.~\cite{Calvetti2000} in 2000.
Many publications on this subject followed; see,~e.g.,~\cite{Calvetti2001,Reichel2005,Baglama2007,Bellalij2015}.

\section{Conclusions}\label{sec:con}

The efforts to develop GMRES algorithms span now three and a half decades.
In this paper we have given an overview of GMRES, including basic algorithms, convergence results, acceleration strategies, parallel algorithms, and a brief account of extensions.
An important area which was not addressed is the performance of GMRES in real-world applications; see,~e.g.,~\cite{Wigton1985,Brown1989,Ohtsuka2004,Liesen2005,Lemieux2008,Gander2015b}, a thorough discussion of which is beyond the scope of this article.

At present, we can see that efficient parallelization of GMRES represents the frontier for linear systems.
Rigorous stability studies of parallel and block algorithms are still lacking; we refer to~\cite{Carson2021} as a recent work in this direction.
Furthermore, for the acceleration strategies mentioned above, preconditioning may offer significant advantages in a variety of contexts, and thus should play a role in solving challenging problems.
There is strong evidence that much effort should be spent in developing mixed-precision parallel GMRES algorithms, along with scalable preconditioners suitable for highly parallel environments.
Describing their finite-precision behavior would be a quite difficult task.

Another area of interest in numerical linear algebra is the design of randomized algorithms; see~\cite{Martinsson2020} for a detailed survey.
We highlight the work by Balabanov and Grigori~\cite{BalabanovX1}, where a randomized Gram-Schmidt algorithm was developed and applied to GMRES.

\bibliography{ref}
\bibliographystyle{abbrv}

\end{document}